\documentclass[12pt]{article}

\usepackage{bussproofs}
\EnableBpAbbreviations

\usepackage{url}

\begin{document}

\title{{\LARGE G\" odel on Deduction}}
\author{{\sc Kosta Do\v sen} and {\sc Milo\v s Ad\v zi\' c}
\\[1ex]
{\small Faculty of Philosophy, University of Belgrade}\\[-.5ex]
{\small \v Cika Ljubina 18-20, 11000 Belgrade, Serbia, and}\\
{\small Mathematical Institute, Serbian Academy of Sciences and Arts}\\[-.5ex]
{\small Knez Mihailova 36, p.f.\ 367, 11001 Belgrade, Serbia}\\[.5ex]
{\small email: kosta@mi.sanu.ac.rs, milos.adzic@gmail.com}}
\date{\small January 2016}
\maketitle

\begin{abstract}
\noindent This is an examination, a commentary, of links between some philosophical views ascribed to G\" odel and general proof theory. In these views deduction is of central concern not only in predicate logic, but in set theory too, understood from an infinitistic ideal perspective. It is inquired whether this centrality of deduction could also be kept in the intensional logic of concepts whose building G\" odel seems to have taken as the main task of logic for the future.
\end{abstract}

\vspace{3ex}

\noindent {\small \emph{Keywords:} deduction, sequent, set, extension, concept, intension, general proof theory, proof-theoretic semantics}

\vspace{3ex}

\noindent {\small \emph{Mathematics Subject Classification
(2010):} 03A05 (Philosophical and critical), 03F03 (Proof theory, general)}

\vspace{3ex}

\noindent {\small \emph{Acknowledgements.} Work on this paper was
supported by the Ministry of Education, Science and Technological Development of Serbia,
while the Alexander von Humboldt Foundation has supported the
presentation of matters related to it by the first-mentioned of us at the conference {\it General Proof Theory:
Celebrating 50 Years of Dag Prawitz's ``Natural Deduction''}, in
T\" ubingen, in November 2015. We are grateful to the organizers for the invitation to the
conference, for their care, and in particular to Peter Schroeder-Heister, for his exquisite hospitality.
We are also grateful to Gabriella Crocco for making G\" odel's
deciphered unpublished notes {\it Max Phil X} available to us, and for allowing us to quote a sentence from them translated in \cite{Ko15}. The publishing of {\it Max Phil X} is part of the project {\it Kurt G\" odel Philosopher: From
Logic to Cosmology}, which is directed by Gabriella Crocco and funded by the French
National Research Agency (project ANR-09-BLAN-0313). She and Antonio Piccolomini d'Aragona were also very kind to invite us to the workshop {\it Inferences and Proofs}, in Marseille, in May 2016, where the first-mentioned of us delivered a talk based partly on this paper, and enjoyed their especial hospitality. We are grateful to the Institute for Advanced Study in Princeton for granting us, through the office of its
librarian Mrs Marcia Tucker and its archivist Mr Casey Westerman, the permission to
quote this sentence; we were asked to give credit for that with the following text:
``All works of Kurt Gödel used with permission. Unpublished Copyright
(1934-1978) Institute for Advanced Study. All rights reserved by Institute
for Advanced Study.''}

\vspace{5ex}

\section{Introduction}
This paper is an attempt to find links between G\" odel's philosophical views and general proof theory (and even proof-theoretic semantics; see the end of Section~4 below). This attempt is not based on the published works of G\" odel, but on the record of his views in \cite{W96}, which cannot be accepted without reservation. It is however quite astonishing how some of the views ascribed to G\" odel accord well with the perspective of general proof theory (see Section~3). Deduction is of central concern for G\" odel not only in predicate logic, but, according to \cite{W96}, he adopts this perspective even for set theory, as a kind of ideal infinitistic extrapolation, inspired to a certain extent by Russell's no class theory (see Section~4).

Words ascribed to G\" odel concerning these matters, which we examine and comment upon here, are not many, but they are quite suggestive. It is to be hoped that the unpublished {\it Nachlass} of G\" odel will yield further confirmation of what these words suggest. It can document at least how G\" odel relied on Gentzen's sequents in the second and last elementary logic course he delivered in his life, which we believe is not a slight matter (see Section~2). We consider at the end whether a link between G\" odel's views and general proof theory should be expected also in the logic of concepts, the intensional logic whose building G\" odel seems to have taken as the main task of logic for the future (see Section~5). Although there is much emphasis in words ascribed to G\" odel concerning this task, not many details about its realization are given.

\section{G\" odel's interest in sequents}
G\" odel's name is one of the most often, if not the most often, mentioned in the field of proof theory that has grown out of Hilbert's program and the quest for consistency proofs. It is however far less often mentioned in general proof theory that has grown out of Gentzen's work on the normal form of proofs, which starts with Gentzen's doctoral thesis \cite{G35}, and has been pursued by Dag Prawitz in \cite{P65}. G\" odel's technical results on intuitionism from the early thirties would be mentioned in this other field of proof theory (among them those on the double-negation translation of \cite{Go33}, which overlap with results independently reached by Gentzen in \cite{G33}), but these results are not of the kind central in the field.

It should not be a surprise that G\" odel paid no doubt close attention to Gentzen's consistency proof for formal arithmetic, in both versions. In \cite{Go38} G\" odel discusses the first published version of \cite{G36}, and notes at the end that the mathematical significance of Gentzen's work is ``tats\" achlich au\ss enordentlich gro{\ss} [in fact extraordinarily great]'',  while the version of \cite{G38} is considered in \cite{Go53} (\S 24). Bernays and Gentzen are called ``leading proponents of formalism'' at the same place (and ``leading formalists'' in the Gibbs lecture \cite{Go51}, p.\ 318, fn 27). Georg Kreisel testifies in \cite{K83} (Note 7(d)(i), p.\ 169): ``[G\" odel] often called Gentzen a better logician than himself.'' If we realize that, according to what is known today, G\" odel in general did not make many compliments concerning logicians who were his contemporaries, what is mentioned in this paragraph gains in importance.

It is however not well known that G\" odel had a favourable opinion about Gentzen's presentation of logic with sequents. It is recorded in \cite{Daw97} (p.\ 136) that in the elementary logic course he gave in the spring of 1939 at the University of Notre Dame his aim was ``to give, as far as possible, a complete theory of logical inferences and logically true propositions, and to show how they can be reduced to a certain number of primitive laws'' (the text under quotes here, and in the remainder of this paragraph, is from G\"odel's notes; see Extract~1 in \cite{CN09}, p.\ 78). After that, one can find at the same place that after giving an axiomatization of propositional logic derived from Russell, and proving in particular completeness, he noted that although Russell's axiomatization had become the standard for textbooks, it is ``open to some objections from the aesthetic point of view'' (see \cite{CN09}, p.\ 70). G\" odel pointed out that statements taken as axioms ``should be as simple and evident as possible... [certainly] simpler than the theorems to be proved, whereas in [Russell's] system... the very simple law of identity $p\supset p$'' was a theorem rather than an axiom (the text quoted here from G\" odel's notes is from \cite{Daw97}; a slightly different rendering is in \cite{CN09}, p.\ 70). It is reported at the same place that G\" odel then turned towards Gentzen's sequents to give an alternative formulation that he regarded as more satisfactory.

This paper, and also the present section, were written without our having been acquainted with G\" odel's unpublished notes for his Notre Dame course. In a companion to this paper \cite{DA16a}, written after gaining this acquaintance, we present G\" odel's system of propositional logic based on sequents and, to provide a context for it, we summarize briefly G\" odel's course, which we present with more detail and comments in another companion \cite{AD16}. We also make there comments upon the course and upon G\" odel's system with sequents, and we reexamine in the light of the notes matters considered in this section.

Besides the elementary logic course at Notre Dame, G\" odel gave an elementary logic course in Vienna in 1935, which is briefly described in \cite{Daw05} (p.\ 153). At the same place, one can find another brief description of the notes for the Notre Dame course, and Gentzen's sequents are again mentioned. They are also mentioned in \cite{Daw} (p.\ 8), where there are further details about G\" odel's time at Notre Dame. Three extracts from G\"odel's notes for the Notre Dame course are published in \cite{CN09}. G\" odel's way of presenting logic with sequents is mentioned in the summary and comments preceding the extracts (p.\ 70, to which we referred above), but the part of the notes with this presentation was not chosen for printing.

Could one assume that G\" odel thought that the advantage of Gentzen's presentation of logic in terms of sequents is {\it merely} aesthetic? As if this aesthetic advantage was a secondary matter, not touching upon the real, primary, goal of logic. For G\" odel, the primary goal of logic, as for mathematics in general, could no doubt be described as reaching the truth about the matters it studies, and this is so from his youth on (see \cite{Go03}, pp.\ 447, point 8, and p.\ 450, point 8(b), and \cite{Go03a}, pp.\ 397-398, 401, 403-405; the gist of these letters of G\" odel may also be found printed in \cite{W74}, pp.\ 8-11). The question is whether reaching this primary goal is separable from reaching the secondary one. Can one reach true, correct, mathematics without at the same time reaching mathematical beauty? Perhaps yes, though one could always blame the failure to reach beauty upon the weakness of men---their inability to recognize it, or to find the right or exact formulation to make it manifest. What however seems impossible is to reach beauty in mathematics without having reached truth. Beautiful incorrect mathematics is out of the question. Its beauty would be false, and not real.

Real beauty is for mathematicians inseparable from truth, but, more than that, it is a sign of truth, and moreover of important truth. This is not a peculiarly G\" odelian point of view, and many, if not most, mathematicians certainly have it. G\" odel has however an explicit relevant statement, taken from a very nice section of his unpublished notes, written a few years after the Notre Dame course:

\begin{quote}
The truth is what has the simplest and the most beautiful \nopagebreak{symbolic expression.}
\end{quote}
(The translation of this sentence from \cite{GoX}, p.\ \textbf{[18]}, is from \cite{Ko15}, Section 4.2.) This sentence could say more than that simplicity and beauty are a sufficient condition for truth, and a sign of it. It could say in addition that they are a necessary condition for truth.

G\" odel has probably kept the same opinion until his old age, according to the following, which are his words as recorded in \cite{W96} (Chapter~4, Section 4.4, p.\ 151), and are nice too:
\begin{quote}
4.4.18 \quad In principle, we can know all of mathematics. It is given to us in its
entirety and does not change---unlike the Milky Way. That part of it of which we
have a perfect view seems beautiful, suggesting harmony; that is, that all the parts
fit together although we see fragments of them only. Inductive inference is not
like mathematical reasoning: it is based on equality or uniformity. But mathematics
is applied to the real world and has proved fruitful. This suggests that the
mathematical and the empirical parts are in harmony and the real world is also beautiful.
Otherwise mathematics would be just an ornament and the real world would be
like an ugly body in beautiful clothing.
\end{quote}

In this vein, which leads to Plato, and seems to have been congenial to G\" odel all his life, truth and beauty, although different, would be indissolubly linked. (Goodness is not mentioned here.) One could perhaps say that, although truth and beauty are intensionally different, they are extensionally equivalent.

G\" odel's sentence first quoted above may also suggest something that presumably was not G\" odel's position---namely, that truth amounts to simplicity and beauty, that it reduces to them, that it is defined by them. As a formalist in the philosophy of mathematics could replace truth by utility, so an {\it aestheticist} in this philosophy, who may, but perhaps need not, be related to the formalist, could replace truth by beauty. The converse reduction of beauty to truth seems to amount to the rejection of beauty as a viable subject in the philosophy of mathematics. In fact, much of the philosophy of mathematics shares with this anti-aestheticism the neglect of beauty.

Be that as it may, it seems safe to say that G\" odel was not making a slight compliment to Gentzen. This matter interests us here, because we may tie it with a rather plausible assumption. The assumption is that G\" odel understood sequents as saying something about deductions, and not as a peculiar way to write implications. This flattened understanding of sequents as implications may be found, for example, in \cite{Ch56} (\S 29 and \S 39.11). Gentzen himself in his thesis \cite{G35} (Sections I.2.4 and  III.1.1) recommends such an understanding, and calls it ``inhaltlich'' (translated as ``informal'' in \cite{G69} and ``intuitive'' in \cite{G55}, but ``soderzhatel'ny\u{\i}'' in \cite{G67}). The connection of sequents with natural deduction, from which they stem, and what Gentzen does with them, suggest strongly that they should be understood as being about deductions, as indeed many of those who worked with them in the last decades understood them.

We speak of deduction, and not consequence, because understanding sequents as being about consequence would not be going far away from the flattened reading of sequents as implications. (Section~4 of \cite{D10} raises the question whether deduction, implication and consequence have {\it inhaltlich} the same meaning.) Strictly speaking, Gentzen does not go beyond taking sequents as partaking in a consequence relation, but his deep concentration on normal form, i.e.\ cut elimination, foreshadows another direction of research. This is the taking of the reflexive, symmetric and transitive closure of the reducibility relation involved in reaching the normal form as a significant equivalence relation on derivations. The first answer to the question of identity criteria for deductions---perhaps the main technical question of general proof theory---is obtained by taking deductions as equivalence classes of derivations (see \cite{D03}, \S 2, and \cite{DP04}, Chapter 1, \S 1.3). The question and the answer are not exactly in Gentzen's writings, but he points very much in that direction. He has a good feeling for the subject, and by pursuing his own goals, which besides a proof of consistency may be also beauty, he prepares well the ground for these investigations that in his time belonged to the future, and are still not entirely accomplished and acknowledged.

How well Gentzen was pointing in the direction of general proof theory becomes clear when one compares his approach to consequence based on sequents to Tarski's approach based on the $Cn$ consequence operation (these matters are considered in \cite{D97}, Section~2, where one can also find references to the relevant papers of Tarski, which are translated in \cite{T56}). In principle, the two formalisms are intertranslatable, but with Tarski's it is hard to imagine how one could reach the idea of normal form for deductions, and how one could adapt $Cn$ to study identity of deduction.

To ascribe to G\" odel the ability to recognize to a certain extent in sequents their potential for doing interesting and nice mathematics about deductions does not seem far-fetched. The same ability could be ascribed rather plausibly to Kleene, somebody from G\" odel's milieu in America (see \cite{Kl83}), who published not much later his book \cite{K52}. Much of Chapter XV of that book is devoted to Gentzen's sequents, and Kleene's paper \cite{K52a}, published the same year, is practically about identity of deductions within Gentzen's sequent systems. The same ability could also be ascribed to Bernays, who supervised Gentzen's work on his thesis \cite{G35} and his later work, and with whom G\" odel corresponded extensively for 45 years (see \cite{Go03}, pp.\ 40-313).

Although G\" odel's contribution to logic was not about deduction, in statements he made in philosophical conversations in his old age, which are recorded in \cite{W96}, deduction became important, and we will consider this matter in the next section.

\section{Logic and deduction}
Let us note first that, instead of the term {\it deduction}, in the words of G\" odel recorded in \cite{W96} one finds rather the term {\it inference}. We will here take these two terms as synonymous. In general, we favour the term {\it deduction}, because we believe that it has less psychologistic overtones (see \cite{D10}, Section~3). We believe that a psychologistic interpretation of inference or deduction does not accord with G\" odel's ideas, as we will try to show at the end of the paper in Section~5.

The most explicit mentioning of the role of deduction, i.e.\ inference, in logic that we have found in \cite{W96} is in the following fragments (from Chapter~8, Section 8.1, pp.\ 266-267):
\begin{quote}
8.4.11 \quad The propositional calculus is about language or deals with the original notion of language: truth, falsity, inference. \ldots

8.4.12 \quad One idea is to say that the function of logic is to allow us to draw inferences. If we define logic by formal evidence directly concerning inference for the finite mind, then there is only one natural choice and it is not natural to treat the infinite as a part of logic. The part of formal inference or formal theory for the finite mind incorporates inferences. The completeness proof of predicate logic confirms its adequacy to this conception of logic.  \mbox{\hspace{.1em}} \hfill\ldots

8.4.14 \quad In contrast to set theory, predicate logic is mainly a matter of rules of inference. It is unnatural to use axioms in it. For the infinite mind, the axioms of set theory are also rules of inference. \ldots

8.4.15 \quad Lower functional calculus [predicate logic] consists of rules of inference. It is not natural to use axioms. It is logic for the finite mind.
\ldots\, For the infinite mind, axioms of set theory are also rules of inference.

8.4.16 \quad For the empiricists, the function of logic is to allow us to draw inferences. It is not to state propositions, but to go over from some propositions to some other propositions. For a theoretical thinker, the propositions embodying such inferences (or implications) are also of interest in themselves.
\end{quote}

The omitted fragment 8.4.13 of this series is about probability, which is not our concern here. We have not quoted the remaining fragments in full, except for the last one. The omitted portions are about matters that, as 8.4.13, are not our concern here: Aristotle, the necessity modal operator, the quantifiers ``most'', ``many'' and ``some (in the sense of plurality)'' (what is exactly this last quantifier is not quite clear), and the need for quantifiers to be able to talk about objects. We omitted from 8.4.14 an interesting claim that according to Bernays mathematics is more abstract than logic, in which we find concepts with content, while in mathematics the concepts of group and field are purely formal. Logic, including its intensional part, the logic of concepts, is characterized in 8.4.18$\,$[part I] (see Section~5 below) as {\it the theory of the formal}.

The term {\it formal} is also used in the quotation from 8.4.12 above. Later in the paper, in Section~5, we will return to what is expressed with this term in that quotation, and also to something else omitted from 8.4.14, which has to do with the logic of concepts. The connection of all that with deduction is not clear.

A warning is now in order concerning the authenticity of these words. A similar warning is made in \cite{W96} (in Chapter~4, p.\ 137) by Hao Wang himself, the author of that book, and in \cite{P98} (Section~5, p.\ 20). The words quoted above were not written by G\" odel. They come from the notes of conversations with G\" odel taken by Hao Wang, who in \cite{W96} says on p.\ 131: ``Following each of these `conversations', I wrote up his ideas
in my own words, and we discussed the fragments I had thus produced''. It seems however that {\it each} in the context where it occurs should not be understood as applying to {\it all} the conversations. Hao Wang says on the same page and on p.\ 137 that the conversations in 1975 and 1976 were usually on the phone, and at the later page he also says: ``I sent him written versions of some
parts of our discussions so that he could comment on them in our next
telephone call. Some of my reports of his sayings in this book are,
therefore in a form that he approved of; other parts he read but was not
satisfied with, even though he recognized them as his own statements. (In
these cases he seemed to expect, or hope, that I would come up with a
clearer exposition than his own.) Still other parts were never seen by him.
Therefore, he probably would not have wished to publish much of the
material in the form I `quote', and it is quite possible that there are places
where I am mistaken about what he actually said.''
On p.\ 130 Hao Wang says: ``I have reconstructed these conversations in several versions, based on very incomplete notes, in an effort (only partially
successful) to interpret them and place them in perspective.'' The publication of the notes in G\" odel's {\it Nachlass} that he wrote in shorthand for himself may show how much Hao Wang was successful. Sometimes the several versions he mentions are not clearly consistent with each other. We are however not aware that what we quoted has been contradicted elsewhere by words ascribed to G\" odel.

The most important question for us that may be raised concerning the quotations above is what is meant there by inference. What formalism is envisaged?

We think that the formalism that fits best the quoted words is that of Gentzen's sequents, which as we have seen in the preceding section interested G\" odel since his young age. An interpretation based on natural deduction would also be possible, but would not be essentially different, since sequents should be taken here as providing only a more precise record of deductions.

The theoretical thinker's ``propositions embodying such inferences (or implications)'' from the last quoted fragment 8.4.16 could be understood as sequents. Moreover, this could well be sequents interpreted as being about deduction, and not about consequence (see the last part of the preceding section). The appearance of the parenthetical {\it implications} may perhaps be taken as a misnomer, coming from a misunderstanding, an unfortunate translation, or a dissatisfaction with the term {\it sequent}, due to G\" odel himself. (This term does not occur in \cite{W96}.) Gentzen, and many others following his usage, write an arrow to separate the left-hand side of sequents from the right-hand side, and this arrow is often used for the implication connective by G\" odel (see \cite{Go86} and \cite{Go90}) and others. As we have seen in the preceding section, Gentzen himself at the beginning of his thesis did not understand sequents as ``inhaltlich'' different from implications.

It is very natural to assume that G\" odel would not think about deductions without thinking also about the rules that justify them, and the ``formal evidence directly concerning inference'' in the quotation from 8.4.12 above may be understood as concerning not only formal sequent systems, but also rules of deduction. Rules of deduction are explicitly mentioned in 8.4.14 and 8.4.15. Could one try to adduce this as ``informal'' evidence for the ability that we ascribed to G\" odel towards the end of the preceding section, which is the ability to recognize in sequents their potential for the investigations of general proof theory? Taking serious account of the rules of deduction is an essential matter in the study of identity criteria for deduction, which we mentioned there, and which is pursued as one of the main tasks of general proof theory. The interesting and important thing for us here is not that G\" odel thought about rules of deduction, but that in 8.4.12 he proposes to ``define logic by formal evidence directly concerning inference for the finite mind''.

An interpretation of our quotations in terms of sequents makes sense everywhere, except for one sentence. Let us first note that what we have quoted from 8.4.14 is quite similar to what we have quoted from 8.4.15, although the former fragment is said to be from 1971 and the latter from 1976. This repetition of the same point after several years may be a sign of its importance for G\" odel, and may corroborate its authenticity, or it may have just been imported in one of the quotations either from the past or from the future. The first part of these quotations about predicate logic fits rather well with the interpretation with sequents understood as being about deductions.

Note then that the last sentence of the first quotation is practically the same as the last sentence of the second quotation. (Does this corroborate authenticity even more, or is it an even stronger sign of importation from the past or the future?) We will consider closely in the next section this sentence, which is difficult to understand.

\section{Set theory and deduction}
The last sentence of our quotation from 8.4.14 in the last section is, modulo the disappearance of a definite article, the same as the last sentence of 8.4.15, and this sentence is:
\begin{quote}
For the infinite mind, [the] axioms of set theory are also rules of inference.
\end{quote}

At first glance, this is obscure. Perhaps some critiques of G\" odel could even accuse him here of mysticism, as he was accused concerning his epistemology of mathematics. (An example of a rather intemperate, but not atypical, accusation of this kind is in \cite{Chi90}, Section I.3, p.\ 21.)

A clue for a reasonable interpretation of this sentence is suggested by the following quotation from \cite{W96} (p.\ 268), which is just a few lines below the quotations of the preceding section:
\begin{quote}
8.4.18$\,$[part II] \quad \ldots\, Elementary logic is the logic for finite minds. If you have an infinite mind, you have set theory. For example, set theory for a finite universe of ten thousand elements is part of elementary logic; compare my Russell paper (probably \cite{Go44}).
\end{quote}
The relevant quotation from \cite{Go44} (p.\ 144, in \cite{Go90}, p.\ 134) to which G\" odel, according to Hao Wang, presumably refers, is the following:
\begin{quote}
\ldots\, \cite{R26} took the course of considering our inability to form propositions of infinite length as a ``mere accident'', to be neglected by the logician. This of course solves (or rather cuts through) the difficulties; but [it] is to be noted that, if one disregards the difference between finite and infinite in this respect, there exists a simpler and at the same time more far reaching interpretation of set theory (and therewith of mathematics). Namely, in case of a finite number of individuals, Russell's {\it aper\c cu} that propositions about classes can be interpreted as propositions about their elements becomes literally true, since, e.g., ``${x\in m}$'' is equivalent to
\[
``x = a_1 \vee x = a_2 \vee \ldots \vee x = a_k"
\]
where the $a_i$ are the elements of $m$; and ``there exists a class such that \ldots'' is equivalent to ``there exist individuals $x_1,x_2,\ldots,x_n$ such that \ldots'', \{footnote 36: The $x_i$ may, of course, as always, be partly or wholly identical with each other.\} provided $n$ is the number of individuals in the world and provided we neglect for the moment the null class which would have to be taken care of by an additional clause. Of course, by an iteration of this procedure one can obtain classes of classes, etc., so that the logical system obtained would resemble the theory of simple types except for the circumstance that mixture of types would be possible. Axiomatic set theory appears, then, as an extrapolation of this scheme for the case of infinitely many individuals or an infinite iteration of the process of forming sets.
\end{quote}
Note that in the last quotation G\" odel speaks of {\it classes}, rather than {\it sets}, which are mentioned only at the end, but at the same time he speaks about ``set theory'', mentioning it twice.

We will not try to reconstruct precisely the reduction of set theory for a finite universe to logic according to this sketch of G\" odel, which comes after his discussion of Russell's ``no class theory'', by which it is inspired. We will however make a few comments on this reduction. (G\" odel's set theory for a finite universe should not be confused with finitary set theory, which amounts to arithmetic, and where the individuals are finite sets, but the universe may be infinite; see \cite{K09} and references therein, and also \cite{C66}, Section I.6).

One might expect that formally this reduction is a {\it translation}, i.e.\ a one-one map, from the language of set theory to a purely logical language. It is not clear what is in general a translation in logic. It is presumably not any one-one map. Sometimes it is required that it should, partly or totally, preserve meaning, and this preservation is often thought to be accomplished by an inductive definition of the map, which makes of it a kind of homomorphism. It is however far from clear that any such inductive definition will guarantee preservation of meaning, which is anyway not a clear concept, as G\" odel might say (see the next section below). It is also not clear that the $\in$-elimination translation suggested by G\" odel should be defined inductively.

To suggest the difficulties of an inductive definition of such a translation $\tau$, let $E$ be the axiom of extensionality:
\[
\forall x\forall y(\forall z(z\in x \leftrightarrow z\in y) \rightarrow x=y),
\]
which should hold in every finite universe of sets, as well as in every universe of sets, and consider how we could obtain $\tau(E)$. If $E'$ is $\forall z(z\in x \leftrightarrow z\in y)$, the antecedent of the implication prefixed in $E$ by $\forall x\forall y$, then $\tau(E')$ would presumably be:
\[
\forall z((z = x_1 \vee \ldots \vee z = x_n) \leftrightarrow (z = y_1 \vee \ldots \vee z = y_n)),
\]
but what would then be $\tau(x=y)$? We cannot just leave $\tau(x=y)$ to be $x=y$, because we want $\tau(E)$ to hold in our universe. We must say in $\tau(x=y)$ something about the $x_i$ and the $y_j$. Perhaps $\tau(x=y)$ should be just $\tau(E')$, which would have the advantage of transforming the axiom of extensionality into the logical truth that is the universal closure of $\tau(E')\rightarrow \tau(E')$. This kind of translation cannot however be given in a context-independent manner, and it is not clear how to define it inductively.

G\" odel's translation $\tau$ would presumably be expected to preserve meaning totally. Does it do that if $\tau(E)$ is the universal closure $\tau(E')\rightarrow \tau(E')$? Do $E$ and $\tau(E)$ have the same meaning? The sentence $\tau(E)$ is tautologous, whereas $E$ does not seem so.

If not equality of meaning, will the translation $\tau$ guarantee at least the equivalence mentioned by G\" odel? What is exactly this equivalence? This is not a straightforward matter. G\" odel suggests that $\tau(\forall x\exists m (x\in m))$ should be the formula:
\[
\forall x \exists x_1 \ldots \exists x_n (x = x_1 \vee \ldots \vee x = x_n).
 \]
The last formula, which can be derived from $x=x$, is logically valid, but $\forall x \exists m(x \in m)$
will not hold if the number of
individuals is finite and there is no infinite descending sequence of elements (as should be the case according to G\" odel's remark in the quotation about the resemblance with type theory). So equivalence cannot mean that for every formula $A$ we have that $A$ holds
in a finite universe if and only if $\tau(A)$ holds in this universe.

What does equivalence mean then? Should one state something weaker? Perhaps that for every formula $A$ we have just the implication that if $A$ holds in a finite universe, then $\tau(A)$ holds in this universe?

Instead of expecting a translation result, one could understand G\" odel's appeal to equivalence just as an admonition: ``Leave aside the mentioning of classes, and mention just individuals.'' Russell's no class theory, although called a {\it theory} (G\" odel seems to call it also an {\it aper\c cu}), could perhaps be understood in a similar manner. To justify this admonition one could rely on something like the implication of the preceding paragraph. Such an implication is sufficient to guarantee that by translating we cannot err, although identity of meaning is doubtful. We cannot fall from truth to falsehood, which is the most important thing for us if we are extensionally-minded.

Anyway, G\" odel in the last two quotations seems to suggest that as set theory for a finite universe reduces to logic, so set theory in general by ``an extrapolation of this scheme for the case of infinitely many individuals'' reduces to logic. And logic is about deduction, as stated in the quotations from the preceding section. The axioms for $\in$ will then be replaced by deduction rules for disjunction, maybe infinite. These are rules like introduction of disjunction and elimination of disjunction in natural deduction, or rules like introduction of disjunction on the right-hand side of a sequent and introduction of disjunction on the left-hand side of a sequent. This is what $\in$ reduces to for an infinite mind.

This is how the sentence quoted at the beginning of this section can be understood reasonably, in a manner pointing towards general proof theory. As the rest of the quotations considered in the preceding section, it could also be taken as pointing in fact towards proof-theoretic semantics (see \cite{SH16}).

Whether the reduction of set theory to logic preserves meaning is however doubtful, as we suggested above. One could perhaps expect an answer in the logic of concepts that G\" odel was hoping for, which we will briefly examine in the next section.

\section{Concepts and deduction}
A great deal of the records of Hao Wang in \cite{W96} of his conversations with G\" odel is about views ascribed to G\" odel concerning the logic that should be developed in the future. This applies in particular to Chapter~8, from which we have been quoting in the preceding two sections, but it applies as well to many other places in the book. On can hope that similar views may be found, at least to a certain extent, in G\" odel's still unpublished, rather considerable, {\it Nachlass}. It would be interesting to systematize this matter, but before the complete publication of the {\it Nachlass}, it might be premature to attempt that. Anyway, it would be too ambitious for us, taking account of our knowledge, and of the space left to us here. All that we venture to do in this closing section is to give a few remarks about the role deduction might play in G\" odel's logic of concepts, the building of which should be the main task of logic for the future. (In \cite{Cro06} G\" odel's views on concepts are considered in the historical perspective provided by related views of Frege and Russell.)

G\" odel was hoping for a logic that will include an intensional theory dealing with concepts, with which he will be able to go beyond the limits of the existing extensional theory. This intensional theory, which will presumably strive for completeness, will not shun away from self-reference, as the extensional theory must do, but it is not clear how the intensional theory will prevent self-referential paradoxes to produce contradiction in it (see \cite{W96}, 8.4.19, p.\ 268, and pp.\ 278-279).

G\" odel found self-reference very productive and took it as a very important tool, with which a theory gains strength. This is suggested by the following quotation (\cite{W96}, p.\ 139):
\begin{quote}
4.3.5 \quad Category theory is built up for the purpose of proving set theory
inadequate. It is more interested in feasible formulations of certain mathematical
arguments which apparently use self-reference. Set theory approaches contradictions to
get its strength.
\end{quote}
This is one of the very rare remarks ascribed to G\" odel that mention category theory. Can one infer from it that category theory is suggested as a guide for the intensional theory of concepts?

We should not relinquish the precious tool that self-reference is if contradiction does not arise.
This attitude of G\" odel is not surprising, since both Cantor's proof of the nonenumerability of the reals, as well as G\" odel's own incompleteness results, involve self-reference in their core. Besides in these two, most famous and important, cases one can no doubt find self-reference applied fruitfully in many other cases. For example, passing from an ordinal $\alpha$ to its successor $\alpha\cup\{\alpha\}$ may be described picturesquely as follows: ``The ordinal $\alpha$ becomes conscious of itself and adds itself to its elements.'' (The finitary set theory of \cite{K09}, mentioned in the preceding section in a parenthetical remark after the quotation from \cite{Go44}, is based on a generalization of this successor operation.)

To accommodate self-reference, G\" odel took Church's formal system of \cite{C32} and \cite{C33}, based on lambda abstraction and application with no types, as a step in the right direction in formulating the logic of concepts (see the references concerning paradoxes before the quotation above). Although Church's system was soon to be shown inconsistent in \cite{KR35}, the possibility remained that one could devise a more satisfactory, consistent, system for the logic of concepts based on similar ideas (see \cite{W96}, 8.4.19, p.\ 268, which we mentioned above).

In the draft of a letter of 1963 written presumably to the theologian Paul Tillich, G\" odel mentions self-reference in the guise of self-knowledge in an even more philosophical context (see \cite{vA06}, pp.\ 259-260; we correct in our quotation some obvious lapses of the draft, and omit most of the editorial apparatus, including footnotes):
\begin{quote}
It occurred to me that in our conversation of last Sunday I answered one of your questions incompletely. I said that in mathematical reasoning the non-computational (i.e. intuitive) element consists in intuitions of higher and higher infinities. This is quite true, but this situation can be further analysed, and then it turns out that they result (as becomes perfectly clear when these things are carried out in detail) from a deeper and deeper self-knowledge of reason [to be more precise from a more and more complete rational knowledge of the {\it essence} of reason (of which essence the faculty of self-knowledge is itself a constituent part)] \ldots
\end{quote}

Let us turn now to the existing extensional theory of logic. It is understood very widely, as in the following quotation, giving the missing part of 8.4.18, of which we have mentioned the later [part II] in the preceding section (\cite{W96}, p.\ 268):
\begin{quote}
8.4.18$\,$[part I] \quad Logic is the theory of the formal. It consists of set theory and the theory
of concepts. The distinction between elementary (or predicate) logic, nonelementary
logic, and set theory is a subjective distinction. Subjective distinctions are
dependent on particular forms of the mind. What is formal has nothing to do with
the mind. Hence, what logic is is an objective issue. Objective logical implication
is categorical. \ldots
\end{quote}

We will first comment upon the first two sentences of this quotation. The extensional theory includes predicate logic and set theory. To treat matters extensionally is the way mathematics does it (\cite{W96}, Section 8.6, p.\ 274):
\begin{quote}
8.6.2 \quad Mathematicians are primarily interested in extensions and we have a
systematic study of extensions in set theory, which remains a mathematical subject
except in its foundations. Mathematicians form and use concepts, but they do not
investigate generally how concepts are formed, as is to be done in logic. We do
not have an equally well-developed theory of concepts comparable to set theory.
At least at the present stage of development, a theory of concepts does not
promise to be a mathematical subject as much as set theory is one.
\end{quote}

The existing extensional predicate logic is about deduction, as we saw in Section~3, and, from an infinitary point of view, set theory is also about deduction, as we saw in the preceding section. So extensionality accords well with deduction. What about the intensional logic of concepts? Will it involve deduction too? The following quotation (\cite{W96}, p.\ 267), which is in the missing part of fragment 8.4.14, from which we quoted in Section~3, makes this doubtful:
\begin{quote}
The whole of set theory is within the purely
formal domain. We have a distinction of two kinds of higher functional calculus
[higher-order logic]: in terms of inferences and in terms of concepts.
\end{quote}
That set theory is formal should mean that it is all within logic, as told in the first two sentences of 8.4.18$\,$[part I] above. Set theory is extensional, and besides it we should have in logic the intensional theory of concepts. If ``in terms of concepts'' stands for ``intensional'', then it is as if ``in terms of inference'' stands for ``extensional''. ``In terms of concepts'' does not then announce an involvement with deduction.

On the other hand, the end of 8.4.18$\,$[part I] above may suggest something else:
\begin{quote}
\ldots\, what logic is is an objective issue. Objective logical implication
is categorical.
\end{quote}
The {\it objective logical implication} mentioned here could perhaps be understood as stemming out of deduction. We have mentioned {\it implication} as a possible misnomer for {\it sequent} in Section~3, in connection with 8.4.16.

Would G\" odel insist here on the {\it objectivity} of implication if he had in mind material implication? The truth value of such an implication is a function of the truth values of the antecedent and the consequent, and is determined out of these truth values quite independently of the human mind. What would it mean to say that material implication is subjective? Who ever thought it is that?

Does G\" odel think here about intuitionistic implication? Maybe he does, but then he might as well have spoken about deduction, to which this implication is closely tied. Intuitionistic implication is characterized through Gentzen's natural-deduction introduction and elimination rules, or sequent rules which amount to these. In terms of category theory, there is behind this characterization an adjoint situation revealed by Lawvere in \cite{Law69}, i.e.\ a natural one-to-one correspondence between deductions from $A\wedge C$ to $B$ and those from $C$ to $A\rightarrow B$, which shows in what precise sense intuitionistic implication is a mirror of deduction (for related matters and further references see \cite{D01} or \cite{D06}).

It would be very surprising if ``categorical'' in the quotation above had to do with this characterization in terms of category theory, although such an interpretation makes perfect sense. It would make more sense, and would be more interesting, than taking this adjective to mean simply something like ``unconditional''. (The characterization in terms of category theory of ``objective'' mathematical structures that deductions make is examined in \cite{D16}.)

The deductions of classical logic associated to Gentzen's plural sequents, with possibly multiple right-hand sides, are studied in terms of category theory in \cite{DP04} (in particular in Chapters 11 and 14). These deductions, which one might sensibly call objective, need not be taken as mirrored by material implication in an adjoint situation, as intuitionistic deductions are mirrored by intuitionistic implication.

Understanding the objectivity of implication, or deduction, as rejecting its psychologistic interpretation makes the quotation above sensible and interesting. This understanding accords well with what precedes this statement of objectivity: ``Subjective distinctions are
dependent on particular forms of the mind. What is formal has nothing to do with
the mind.'' The rejection of psychologism with respect to deduction is something we mentioned briefly at the beginning of Section~3 (and we gave a reference concerning that matter).

Will objective deductions be represented in the logic of concepts? Since this logic strives for a complete characterization of concepts, it should not leave out the objective deductions tied to concepts. Will these deductions be derived from another characterization of concepts, or must we concentrate on deductions directly?

One might suppose that a possible clue for connecting concepts with deductions is provided by G\" odel's views on the concept of absolute proof. A particularly telling remark concerning this matter is the following (\cite{W96}, p.\ 188):
\begin{quote}
6.1.13 \quad The concept of {\it concept} and the concept of {\it absolute proof} may be mutually definable. \ldots
\end{quote}
It is however not clear that the notion of absolute proof is tied specifically to deduction. It seems rather to be something that should supersede provability in a formal system, i.e.\ recursive enumerability. This is corroborated by G\" odel's discussion of absolute provability in the introductory part of \cite{Go46}.

The concepts of extensional logic, concepts tied to the connectives, the quantifiers and equality, involve deduction, according to what was said in Section~3. According to what was said in Section~4, the same applies, in the perspective of an infinite mind, to the concept {\it set}, another concept of extensional logic. Shouldn't then that part of the logic of concepts that deals with these concepts, the intensional logic of these concepts, if not of others, involve deduction?

\end{document}